\documentclass[11pt]{article}

\usepackage{rotating}
\usepackage{etex}
\reserveinserts{28}
\usepackage[all]{xy}
\usepackage{enumitem}

\usepackage{amsfonts}
\usepackage{amsthm}
\usepackage{graphicx}
\usepackage{mathrsfs}
\usepackage{bm}
\usepackage{cite}
\usepackage[colorlinks,linkcolor=blue,citecolor=blue]{hyperref}
\usepackage{amssymb,amsmath} % font used for R in Real numbers
\usepackage{hyperref}
\usepackage{makecell}
\usepackage{pgf}
\usepackage{pgffor}
\usepackage{pgfcalendar}
\usepackage{pgfpages}
\usepackage{shuffle,yfonts}
\usepackage{mathtools}
\usepackage{tikz}
\usetikzlibrary{arrows,shapes,chains,patterns}
\allowdisplaybreaks

\newcommand{\FES}{\mathsf {FES}}
\newcommand{\FMZV}{\mathsf {FMZV}}
\newcommand{\FMtV}{\mathsf {FMtV}}
\newcommand{\FMTV}{\mathsf {FMTV}}
\newcommand{\FMSV}{\mathsf {FMSV}}
\newcommand{\FMMV}{\mathsf {FMMV}}
\newcommand{\FMMVo}{\mathsf {FMMVo}}
\newcommand{\FMMVe}{\mathsf {FMMVe}}
\newcommand{\FAMMV}{\mathsf {FAMMV}}

\newcommand{\ES}{\mathsf {ES}}
\newcommand{\MZV}{\mathsf {MZV}}
\newcommand{\MtV}{\mathsf {MtV}}
\newcommand{\MTV}{\mathsf {MTV}}
\newcommand{\MSV}{\mathsf {MSV}}
\newcommand{\MMV}{\mathsf {MMV}}
\newcommand{\MMVo}{\mathsf {MMVo}}
\newcommand{\MMVe}{\mathsf {MMVe}}

\newcommand{\sha}{\shuffle}

\newcommand{\Sy}{{\mathcal S}}

\newcommand{\gb}{\beta}
\newcommand{\gc}{\gamma}
\newcommand{\gd}{\delta}

\newcommand{\ola}{\overleftarrow}

\newcommand\frakS{{\mathfrak{S}}}

\newcommand\mya{{\texttt{a}}}
\newcommand\myb{{\texttt{b}}}
\newcommand\myc{{\texttt{c}}}

\newcommand{\tq}{{\texttt{q}}}
\newcommand{\tG}{{G}}

\newcommand{\emtq}{{\emph{\texttt{q}}}}

\newcommand\eps{{\varepsilon}}

\newcommand{\bfs}{{\boldsymbol{\sl{s}}}}

\newcommand\bfsi{{\boldsymbol \sigma}}

\newcommand\bfeps{{\boldsymbol \varepsilon}}
\newcommand\bfgs{{\boldsymbol \sigma}}

\newcommand{\calA}{\mathcal{A}}
\newcommand{\calB}{\mathcal{B}}
\newcommand{\calM}{\mathcal{M}}
\newcommand{\calP}{\mathcal{P}}
\newcommand{\calT}{\mathcal{T}}

\catcode`!=11
\let\!int\int \def\int{\displaystyle\!int}
\let\!lim\lim \def\lim{\displaystyle\!lim}
\let\!sum\sum \def\sum{\displaystyle\!sum}
\let\!sup\sup \def\sup{\displaystyle\!sup}
\let\!inf\inf \def\inf{\displaystyle\!inf}
\let\!cap\cap \def\cap{\displaystyle\!cap}
\let\!max\max \def\max{\displaystyle\!max}
\let\!min\min \def\min{\displaystyle\!min}
\let\!frac\frac \def\frac{\displaystyle\!frac}
\catcode`!=12

\let\oldsection\section
\renewcommand\section{\setcounter{equation}{0}\oldsection}

\allowdisplaybreaks

\DeclareMathOperator*{\sgn}{sgn}
\DeclareMathOperator*{\dep}{dep}

\DeclareMathOperator{\Coeff}{Coeff}

\def\gs{{\sigma}}

\def\N{\mathbb{N}}
\def\Z{\mathbb{Z}}
\def\Q{\mathbb{Q}}

\def\db{\mathbb{D}}

\def\ol{\overline}

\theoremstyle{plain}
\newtheorem{thm}{Theorem}[section]

\newtheorem{conj}[thm]{Conjecture}

\newtheorem{prop}[thm]{Proposition}
\theoremstyle{definition}
\newtheorem{defn}{Definition}[section]
\newtheorem{rem}[thm]{Remark}

\newtheorem{prob}[thm]{Problem}

\begin{document}
%%%%%%%%%%%%%%%%%%%% title %%%%%%%%%%%%%%%%%%%%%%%%%%%%%%%%%%%%%%%%%%%%%%%%
\title{\bf Finite Multiple Mixed Values} 

\author{Jianqiang Zhao\thanks{Email: zhaoj@ihes.fr}\\[1mm]
\ \\
with an Appendix by J.\ Feng, A.\ Kim, S.\ Li, R.\ Qin,  L.\ Wang, and J.\ Zhao\\
\ \\
\small   Department of Mathematics, The Bishop's School, La Jolla, CA 92037, USA}

\date{}
\maketitle

\noindent{\bf Abstract.} In recent years, a variety of variants of multiple zeta values (MZVs) have been defined
and studied. One way to produce these variants is to restrict the indices in the definition of MZVs to some fixed
parity pattern, which include Hoffman's multiple $t$-values, Kaneko and Tsumura's multiple $T$-values, and
Xu and the author's multiple $S$-values. We have also considered the so-called multiple mixed values by allowing 
all possible parity patterns and studied a few important relations among these values. In this paper, we will 
turn to their finite analogs and their symmetric forms, motivated by a deep conjecture of Kaneko and Zagier which
relates the finite MZVs and symmetric MZVs, and a generalized version of this conjecture by the author to the Euler
sum (i.e., level two) setting. We will present a few important relations among these values such as the
stuffle, reversal, and linear shuffle relations. We will also compute explicitly the (conjecturally smallest) generating set in weight one and two cases. In the appendix we tabulate some dimension computation for various sub-spaces of the finite multiple mixed values and propose a conjecture.

\medskip



\section{Introduction}
\subsection{Multiple zeta values and their finite analogs}
For any composition of positive integers $\bfs=(s_1,\dots,s_d)\in\N^d$, we define the \emph{multiple zeta value} (MZV) by
\begin{equation}\label{equ:psumMZV}
\zeta(\bfs):=\sum_{n_1>\dots>n_d>0} \prod_{j=1}^d \frac{1}{n_j^{s_j}}
\end{equation}
and the \emph{multiple zeta star value} (MZSV) by
\begin{equation}\label{equ:psumMZSV}
\zeta^\star(\bfs):=\sum_{n_1\ge \dots\ge n_d\ge 1}  \prod_{j=1}^d \frac{1}{n_j^{s_j}}.
\end{equation}
These converge if and only if $s_1\ge 2$ in which case we say $\bfs$ is \emph{admissible}. As usual, we call
$|\bfs|:=s_1+\cdots+s_d$ the weight and $d$ the depth. These values were first systematically studied
by Zagier \cite{Zagier1994} and Hoffman \cite{Hoffman1992} independently.

In recent years, a lot of research has been done concerning the structure of different variants of multiple zeta values
due to their important applications in both mathematics and theoretical physics. For example, when we allow alternating
signs to appear then we obtain the so-called Euler sums (also called alternating MZVs) which play important roles
in the study of knot theory and  Witten multiple zeta function associated with Lie algebra (see, e.g., \cite{Broadhurst1996,Zhao2011b}).
If we allow not only $\pm 1$ but more generally $N$th roots of unity then we can consider colored/cyclotomic MZVs of level $N$,
which have also appeared unexpectedly in the study of Feynman diagrams \cite{Broadhurst1999}.

On the other hand, the modular arithmetic nature of the partial sums of MZVs was first considered by Hoffman \cite{Hoffman2015} and the last author \cite{Zhao2008a} independently. Contrary to the classical cases above, not many variants of these sums exist. To set up the correct
theoretical framework for these variants, we define the following ad\'ele-like ring which was first considered by Kontsevich and then applied
to the $p$-adic setting by Kaneko and Zagier \cite{KanekoZa2013}. Let $\calP$ be the set of primes. Set
\begin{equation}\label{equ:poormanA}
\calA:=\prod_{p\in\calP}(\Z/p\Z)\bigg/\bigoplus_{p\in\calP}(\Z/p\Z).
\end{equation}
Then we can define the \emph{finite multiple zeta values} (FMZVs) by the following:
\begin{equation}\label{defn:FMZV}
\zeta_\calA(\bfs):=\left(H_p(\bfs):=\sum_{p>n_1>\dots>n_d>0}
 \prod_{j=1}^d \frac{1}{n_j^{s_j}} \pmod{p}\right)_{p\in\calP}\in \calA.
\end{equation}
In 2014, Kaneko and Zagier proposed a deep conjecture (see Conjecture~\ref{conj:KanekoZagierAltVersion} below for a generalization)
relating these values on the $p$-adic side
to MZVs on the Archimedean side via a mysterious connection. This conjecture is far from being proved but since then
a plethora of parallel results have been shown to hold on both sides simultaneously (see, e.g., \cite{Murahara2016,SaitoWa2013b,Sakurada2023,SingerZ2019}).
In particular, for each positive integer $w\ge 2$, the element
\begin{equation}\label{equ:betaw}
   \gb_w:=\Big(\frac{B_{p-w}}{w}\Big)_{w<p\in\calP} \in \calA
\end{equation}
is the finite analog of $\zeta(w)$, where $B_n$'s are the Bernoulli numbers defined by
\begin{equation*}
    \frac{t}{e^t-1}=\sum_{n\ge 0}B_n\frac{ t^n}{n!}.
\end{equation*}
And the Fermat quotient
\begin{equation}\label{equ:FermatQ2}
  \tq_2:=\Big(\frac{2^{p-1}-1}{p}\Big)_{2<p\in\calP} \in \calA
\end{equation}
is the analog of $-\zeta(\bar1)=\log 2$.  

\subsection{Euler sums and their finite analogs}
For $s_1,\dots,s_d\in\N$ and $\eps_1,\dots,\eps_d=\pm 1$, we define the \emph{Euler sums}
\begin{equation}\label{equ:z}
\zeta\binom{s_1,\dots,s_d}{\eps_1,\dots,\eps_d}
   := \sum_{n_1>\cdots>n_d>0}\;\prod_{j=1}^d  \frac{\eps_j^{n_j}}{n_j^{s_j}}.
\end{equation}
To save space, if $\eps_j=-1$ then $\overline s_j$ will be used and if a substring $S$ repeats $n$
times in the list then $\{S\}^n$ will be used. For example, $\zeta(\{3\}^n)=8^n\zeta(\{\bar2,1\}^n)$ for all $n\in\N$
(see \cite{Zhao2010a}).

To state the relations between Euler sums concisely we define a kind of double cover of the set $\N$ of positive integers.
\begin{defn} \label{defn:dbANDoplus}
Let $\db$ be the set of \emph{signed numbers}
$\N \cup \ol{\N}$ where
\begin{equation*}
  \ol{\N}:=\{\bar s: s\in\N\}.
\end{equation*}
Define the absolute value function $| \cdot |$ on $\db$ by $|s|=|\bar s|=s$
for all $s\in\N$ and the sign function by $\sgn(s)=1$ and $\sgn(\bar s)=-1$
for all $s\in\N$.  We make $\db$ a semi-group by defining
a commutative and associative
binary operation $\oplus$ (called \emph{O-plus})
as follows: for all $a,b\in\db_0$
\begin{equation}\label{equ:oplusDefn}
  a\oplus b:=
\left\{
 \begin{array}{ll}
  \ol{|a|+|b|}, \quad&\quad \hbox{if $\sgn(a)\ne \sgn(b)$;} \\
  |a|+|b|, \quad & \quad\hbox{if $\sgn(a)=\sgn(b)$.}
 \end{array}
\right.
\end{equation}
Finally, we define the conjugate operation $\ol{s}$ toggling between $\N$ and $\ol{\N}$.
\end{defn}

For any $\bfs=(s_1,\ldots, s_d)\in\db^d$, we define the $n$-th partial sum of the \emph{Euler sums} by
\begin{equation}\label{equ:psumES}
\zeta_n(\bfs):=\sum_{n>n_1>\dots>n_d>0} \prod_{j=1}^d \frac{\sgn(s_j)^{n_j}}{n_j^{|s_j|}} .
\end{equation}
Similarly to FMZVs, finite \emph{Euler sums} (FESs) are defined by
\begin{equation}\label{defn:FES}
\zeta_\calA(\bfs):=\Big(\zeta_p(\bfs) \pmod{p}\Big)_{p\in\calP}\in \calA.
\end{equation}

In \cite[Conjecture~8.6.9]{ZhaoBook} the author of this paper extended Kakeko--Zagier conjecture to the setting of
the Euler sums. For $\bfs=(s_1,\ldots, s_d)\in\db^d$, define the symmetrized version of the alternating Euler sums by
\begin{align*}
\zeta_\ast^\Sy(\bfs):=&\sum_{i=0}^d
\left(\prod_{j=1}^i (-1)^{|s_j|} \sgn (s_j)\right) \zeta_\ast(s_i,\dots,s_1) \zeta_\ast(s_{i+1},\dots,s_d),\\
\zeta_\sha^\Sy(\bfs):=&\sum_{i=0}^d
\left(\prod_{j=1}^i (-1)^{|s_j|} \sgn (s_j)\right) \zeta_\sha(s_i,\dots,s_1) \zeta_\sha(s_{i+1},\dots,s_d),
\end{align*}
where $\zeta_\sharp$ ($\sharp=\ast$ or $\sha$) are regularized values (see \cite[Proposition~13.3.8]{ZhaoBook}).
They are called $\sharp$-\emph{symmetric Euler sums}. If $\bfs\in\N^d$ then they are called
$\sharp$-\emph{symmetric multiple zeta values} ($\sharp$-SMZVs or simply SMZVs if $\sharp$ does not matter).

\begin{conj}\label{conj:KanekoZagierAltVersion}
For any $w\in\N$, let $\FES_{w}$ (resp.\ $\ES_w$) be the $\Q$-vector space generated by
all FESs (resp.\ Euler sums) of weight $w$. Then there is an isomorphism
\begin{align*}
f_\ES: \FES_{w} & \longrightarrow \frac{\ES_w}{\zeta(2)\ES_{w-2}} , \\
 \zeta_\calA(\bfs) & \longmapsto \zeta_\sha^\Sy(\bfs).
\end{align*}
\end{conj}

\begin{rem}
Note that we can replace $\zeta_\sha^\Sy(\bfs)$ by $\zeta_\ast^\Sy(\bfs)$ (see \cite[Exercise 8.7]{ZhaoBook}).
\end{rem}

To better understand this mysterious relation is the primary motivation of this paper. We will mainly study a few variants of
finite analogs of Euler sums by presenting some results that are analogous to
those on the Archimedean side.

\section{Multiple mixed values and their finite analogs}
Put $\Z_*=\Z_{\ne 0}$. For any $\bfs\in \Z_*^d$ and $n\in\N$ we define the $n$th partial sum of \emph{multiple mixed values} (MMVs) by
\begin{equation}\label{defn:psumMMV}
 M_n(\bfs):=\sum_{n>n_1>\dots>n_d>0} \prod_{j=1}^d \frac{1+\sgn(s_j)(-1)^{n_j}}{2n_j^{|s_j|}}
=\sum_{\substack{n>n_1>\dots>n_d>0\\ n_j\equiv (1-\sgn(s_j))/2 \pmod{2} }} \prod_{j=1}^d \frac{1}{n_j^{|s_j|}}
\end{equation}
and the \emph{finite multiple mixed values} (FMMVs) by
\begin{equation}\label{defn:FMMV}
 M_\calA(\bfs):=\Big(  M_p(\bfs) \pmod{p}\Big)_{p\in\calP}\in\calA.
\end{equation}
We call $|\bfs|:=|s_1|+\cdots+|s_d|$ the weight and $d$ the depth.

The motivation to define the MMVs is to find a common generalization of a few variants of level two MZVs
including the following. For all admissible $\bfs=(s_1,\dots,s_d)\in\N^d$, we define the
\emph{\emph{multiple $t$-values}} (MtVs, see \cite{Hoffman2019}),
\emph{multiple $T$-values} (MTVs, see \cite{KanekoTs2020}),
and \emph{multiple $S$-values} (MSVs, see \cite{XuZhao2020a}) by
\begin{align}\label{defn:MtV}
 t(\bfs):=&\, M(-s_1,\dots,-s_j,\dots,-s_d)=\sum_{n_1>\dots>n_d>0} \prod_{j=1}^d \frac{1}{(2n_j-1)^{s_j}},\\
\label{defn:MTV}
T(\bfs):=&\, M\Big((-1)^d s_1,\dots,(-1)^{d-j+1}s_j,\dots,-s_d\Big)=\sum_{\substack{n_1>\dots>n_d>0\\ n_j\equiv d-j+1\pmod{2}}} \prod_{j=1}^d \frac{1}{n_j^{s_j}},\\
\label{defn:MSV}
S(\bfs):=&\, M\Big((-1)^{d-1} s_1,\dots,(-1)^{d-j}s_j,\dots,s_d\Big)=\sum_{\substack{n_1>\dots>n_d>0\\ n_j\equiv d-j\pmod{2}}} \prod_{j=1}^d \frac{1}{n_j^{s_j}},
\end{align}
respectively. Their finite analogs  $t_\calA(\bfs)$, $T_\calA(\bfs)$, and $S_\calA(\bfs)$
are defined similarly as in \eqref{defn:FMZV} and \eqref{defn:FES}. It is clear that
\begin{equation*}
t_\calA(\bfs)=\frac{1}{2^d}\sum_{\eps_1,\dots,\eps_d=\pm 1}\bigg(\prod_{1\le j\le d} \eps_j\bigg)  \zeta_\calA\binom{\bfs}{\bfeps},\quad
F_\calA(\bfs)=\frac{1}{2^d}\sum_{\eps_1,\dots,\eps_d=\pm 1}\bigg(\prod_{\substack{1\le j\le d\\ 2|d-j\ \text{if $F=T$}\\  2\nmid d-j\ \text{if $F=S$}}} \eps_j\bigg)  \zeta_\calA\binom{\bfs}{\bfeps}.
\end{equation*}
More recently, another type of level two FMZVs is defined by Kaneko et al.\ \cite{KanekoMuYo2021} as follows
(after multiplying by $2^{|\bfs|}$ on the following):
\begin{align}\label{equ:finitezeta2Def}
\zeta_\calA^{(2)}(\bfs):=\bigg( \sum_{p>n_1>\cdots>n_d>0,\ 2|n_j\, \forall j} \ \prod_{j=1}^d \frac{1}{n_j^{s_j}}\bigg)_{p\in\calP}
= \frac{1}{2^d} \sum_{\eps_1,\dots,\eps_d=\pm 1}  \zeta_\calA\binom{\bfs}{\bfeps}.
\end{align}
It is also clear that 
\begin{align}\label{equ:zeta2Def}
\zeta^{(2)}(\bfs):=\sum_{n_1>\cdots>n_d>0,\ 2|n_j\, \forall j} \ \prod_{j=1}^d \frac{1}{n_j^{s_j}} 
= \frac{1}{2^{|\bfs|}} \zeta(\bfs).
\end{align}

\begin{defn}
Let $d\in\N$ and $\bfs=(s_1,\dots,s_d)\in\N^d$. We define the \emph{$\sharp$-regularized MTVs}
($\sharp=\ast$ or $\sha$) and \emph{MSVs} by
\begin{align*}
F_\sharp(\bfs):=&\, \frac{1}{2^d}\sum_{\gs_1,\dots,\gs_d=\pm 1}\bigg(\prod_{\substack{1\le j\le d\\ 2|d-j\ \text{if $F=T$}\\  2\nmid d-j\ \text{if $F=S$}}} \gs_j\bigg)  \zeta_\sharp(\bfs;\bfgs)  \quad (\text{$F=T$ or $S$}).
\end{align*}
We define the \emph{$\sharp$-symmetric multiple $T$-values} (SMTVs) and
 \emph{$\sharp$-symmetric multiple $S$-values} (SMSVs) by
\begin{equation*}
F_\sharp^\Sy(\bfs):=
\left\{
  \begin{array}{ll}
  \displaystyle \sum_{i=0}^d  \Big(\prod_{\ell=1}^i (-1)^{s_\ell} \Big) F_\sharp(s_i,\dots,s_1)F_\sharp(s_{i+1},\dots,s_d),\quad & \quad \hbox{if $d$ is even;} \\
  \displaystyle \sum_{i=0}^d  \Big(\prod_{\ell=1}^i (-1)^{s_\ell} \Big) \widetilde{F}_\sharp(s_i,\dots,s_1)F_\sharp(s_{i+1},\dots,s_d), \quad&\quad \hbox{if $d$ is odd,}
  \end{array}
\right.
\end{equation*}
where $\widetilde{F}=S+T-F$ and we set $\prod_{\ell=1}^0=1$ as usual.
\end{defn}

In \cite{Zhao2024b} we proved that for $\sharp=\ast$ or $\sha$ and for all $\bfs\in\N^d$
\begin{equation}\label{equ:fESofTS}
f_\ES T_\calA(\bfs)=T_\sharp^\Sy(\bfs) \quad \text{and}\quad f_\ES S_\calA(\bfs)=S_\sharp^\Sy(\bfs) \pmod{\zeta(2)}.
\end{equation}

In comparison to the results in \cite{XuZhao2020a}, we define the following subspaces of $\FES$:
\begin{itemize}
  \item $\FES$: generated by FESs
  \item $\FMMV$: generated by FMMVs
  \item $\FMtV$: generated by FMtVs
  \item $\FMTV$: generated by FMTVs
  \item $\FMSV$: generated by FMSVs
  \item $\FMZV^{(2)}$: generated by level two FMZVs $\zeta_\calA^{(2)}(\bfs)$ defined by \eqref{equ:finitezeta2Def}
  \item $\FMMVe$: generated by the $\{M_\calA(\bfs): \bfs\in\Z_*^d,s_d>0\}$
  \item $\FMMVo$: generated by the $\{M_\calA(\bfs): \bfs\in\Z_*^d,s_d<0\}$
\end{itemize}

\begin{prop}\label{prop:FES=FMMV}
We have $\FES_w=\FMMV_w$ for all $w\in\N$.
\end{prop}
\begin{proof}
This follows easily from the identities
\begin{equation*}
(-1)^n=\frac{(1+(-1)^n)-(1-(-1)^n)}2,\quad 1=\frac{(1+(-1)^n)+(1-(-1)^n)}2.
\end{equation*}
\end{proof}

For any fixed weight $w\ge 1$ we clearly have the inclusion relations between weight $w$ pieces of the above subspaces:
$$\begin{array}{ccccccccc}
   \   & \  & \FMSV_w & \subseteq & \FES_w & \supseteq & \FMTV_w & \  &  \   \\
   \
   &  \    & \text{\raisebox{-2pt}{\begin{turn}{90} $\supseteq$ \end{turn}}}
   &  \    & \text{\raisebox{-2pt}{\begin{turn}{90} $=$ \end{turn}}}
   &  \    & \text{\raisebox{-2pt}{\begin{turn}{90} $\supseteq$ \end{turn}}}
   &  \    & \     \\
\FMZV^{(2)}_w &\subseteq& \FMMVe_w&\subseteq&\FMMV_w &\supseteq& \FMMVo_w&\supseteq& \FMtV_w
\end{array}$$

Kaneko et al. \cite[Conjecture~5]{KanekoMuYo2021} conjectured that $\FMZV^{(2)}_w=\FES_w$ for all $w\in\N$.
We further conjecture the equal signs in the second row above always hold (also see Conjecture~\ref{conj:dimFMMVs}).

\begin{conj} \label{conj:RelBetweenSubSpaces}
For sufficiently large weight $w\in\N$, we have
$$\begin{array}{ccccccccc}
   \FMZV_w    & \subsetneq & \FMSV_w & \subsetneq & \FES_w & \supsetneq & \FMTV_w & \supsetneq &  \FMZV_w  \\
   \text{\raisebox{-2pt}{\begin{turn}{90} $\supsetneq$ \end{turn}}}
   &  \    & \text{\raisebox{-2pt}{\begin{turn}{90} $\supsetneq$ \end{turn}}}
   &  \    & \text{\raisebox{-2pt}{\begin{turn}{90} $=$ \end{turn}}}
   &  \    & \text{\raisebox{-2pt}{\begin{turn}{90} $\supsetneq$ \end{turn}}}
   &  \    & \text{\raisebox{-2pt}{\begin{turn}{90} $\supsetneq$ \end{turn}}}     \\
\FMZV^{(2)}_w &=& \FMMVe_w&=&\FMMV_w &=& \FMMVo_w&=& \FMtV_w
\end{array}$$
\end{conj}

Except for the middle vertical equal sign all the other relations are supported by numerical evidence but no formal proofs yet.
In contrast, we have
the following conjectural relations between the classical subspaces of Euler sums (cf. the Venn diagram at the end of \cite{XuZhao2020a}).

\begin{prob}
How can we verify numerically the inclusion $\FMZV_w  \subseteq  \FMSV_w \cap \FMZV^{(2)}_w$ and $\FMZV_w  \subseteq  \FMTV_w \cap \FMtV_w$?
We only need to find a basis in each of the four subspaces on the right and show that every FMZV can be expressed as a $\Q$-linear
combination of the basis elements. Is it even possible that both (or one) of the inclusions are actually equalities?
\end{prob}

\begin{conj} For sufficiently large weight $w\in\N$, we have
$$\begin{array}{ccccccccc}
   \MZV_w    & \subsetneq & \MSV_w & \subsetneq & \ES_w & \supsetneq & \MTV_w & \supsetneq &  \MZV_w  \\
   \text{\raisebox{-2pt}{\begin{turn}{90} $=$ \end{turn}}}
   &  \    & \text{\raisebox{-2pt}{\begin{turn}{90} $\supsetneq$ \end{turn}}}
   &  \    & \text{\raisebox{-2pt}{\begin{turn}{90} $\subsetneq$ \end{turn}}}
   &  \    & \text{\raisebox{-2pt}{\begin{turn}{90} $\supsetneq$ \end{turn}}}
   &  \    & \text{\raisebox{-2pt}{\begin{turn}{90} $\supsetneq$ \end{turn}}}     \\
\MZV^{(2)}_w &\subsetneq & \MMVe_w&= &\MMV_w &\supsetneq& \MMVo_w&\supsetneq& \MtV_w
\end{array}$$
\end{conj}
The left vertical equal sign is obvious by \eqref{equ:zeta2Def} and the middle vertical (strict) inclusion follows from \cite[Theorem~7.1]{XuZhao2020a}
if we assume a variant of Grothendieck's period conjecture.
Again, all the other relations are supported by numerical evidence but no formal proofs yet.

\subsection{Stuffle relations}
Stuffle relations hold in all of the above subalgebras of $\FES$ in Conjecture~\ref{conj:RelBetweenSubSpaces}
except $\FMTV$ and $\FMSV$. For example, for all $a,b\in\db$
$$\zeta_\calA(a)\zeta_\calA(b)=\zeta_\calA(a,b)+\zeta_\calA(b,a)+\zeta_\calA(a\oplus b),$$
and if $a,b,c\in\Z_*$
\begin{align*}
M_\calA(a,b)M_\calA(c)=&\, M_\calA(a,b,c)+M_\calA(a,c,b)+M_\calA(c,a,b)\\
&\, +\gd_{\sgn(a),\sgn(c)}M_\calA(a+c,b)+\gd_{\sgn(b),\sgn(c)}M_\calA(a,b+c)
\end{align*}
where $\gd_{s,t}$ is the Kronecker symbol satisfying $\gd_{s,t}=1$ if $s=t$ and $\gd_{s,t}=0$ otherwise.
Hence the stuffing for FMMVs occurs only when the two merging components have the same sign.

\begin{prob}
Do $\FMTV$ and $\FMSV$ form subalgebras of $\FES$? If not, what is the first instance of a product that is not closed?
Note, if the answer is negative then it cannot be verified rigorously with current level of knowledge.
This is similar to the situation for classical $\MTV$ and $\MSV$ and it is for
same type of reason. In the classical setting, counterexamples can only be verified numerically approximately but they cannot be proved
rigorously due to difficulty in proving transcendence result. In the $p$-adic cases, the transcendence problem for $\calA$-numbers may be even harder. For example, we don't even know if $\gb_w\in\calA$ vanishes or not for any odd $w\ge 3$ although conjecturally
the density of primes $p>w$ such that $B_{p-w}\equiv 0 \pmod{p}$ should be 0.
\end{prob}

Here is the process to find such a product of two FMTVs of total weight $w$ numerically, under the assumption that $\dim_\Q \FES_w=F_w$ which is the same as $\dim_\Q \FMMV_w=F_w$ by Proposition~\ref{prop:FES=FMMV}. First, find a generating set $\calT\calB_w$ of $\FMTV_w$ which is possible by using linear shuffle, reversal and stuffle relations. Second, expanding the set into a a generating set $\calM\calB_w$ of $\FMMV_w$. If $|\calM\calB_w|=F_w$ then it is a basis by the assumption $\dim_\Q \FMMV_w=F_w$. Finally,
for each product of FMTVs of weight $w$ we can express it using the conjectural basis $\calM\calB_w$. If
elements outside of $\calT\calB_w$ are needed for such a product then it does not lie in $\FMTV_w$,
meaning the product of $\FMTV$ is not closed.

\subsection{Reversal Relations}

\begin{prop} \label{prop:FMMVreversal}
For all $d\in\N$ and $\bfs\in\Z_*^d$ then
\begin{align}\label{equ:FMMVReversal}
M_\calA(\ola{\bfs})=&\, (-1)^{|\bfs|} M_\calA(-\bfs).
\end{align}
In particular, for all $\bfs\in\N^d$
\begin{align}\label{equ:FMtSVReversalEvenDep}
t_\calA(\ola{\bfs})=&\, (-1)^{|\bfs|} \zeta^{(2)}_\calA(\bfs), \\
T_\calA(\ola{\bfs})=&\, (-1)^{|\bfs|} T_\calA(\bfs) \quad\text{and}\quad S_\calA(\ola{\bfs})=(-1)^{|\bfs|} S_\calA(\bfs) \quad\text{if }  2|d, \label{equ:FMTSVReversalEvenDep} \\
T_\calA(\ola{\bfs})=&\, (-1)^{|\bfs|} S_\calA(\bfs)\quad\text{and}\quad S_\calA(\ola{\bfs})=(-1)^{|\bfs|} T_\calA(\bfs) \quad\text{if } 2\nmid d.\label{equ:FMTSVReversalOddDep}
\end{align}
\end{prop}

\begin{proof}
The relations follow immediately from the change of indices $n\to p-n$:
\end{proof}

\begin{prop} \label{prop:SMMVreversal}
For all $d\in\N$ and $\bfs\in\N^d$
\begin{align}\label{equ:SMTSVReversalEvenDep}
T_\ast^\Sy(\ola{\bfs})=&\, (-1)^{|\bfs|} T_\ast^\Sy(\bfs) \quad\text{and}
\quad S_\ast^\Sy(\ola{\bfs})=(-1)^{|\bfs|} S_\ast^\Sy(\bfs) \quad\text{if } 2|d, \\
T_\ast^\Sy(\ola{\bfs})=&\, (-1)^{|\bfs|} S_\ast^\Sy(\bfs)\quad\text{and}
\quad S_\ast^\Sy(\ola{\bfs})=(-1)^{|\bfs|} T_\ast^\Sy(\bfs) \quad\text{if } 2\nmid d. \label{equ:SMTSVReversalOddDep}
\end{align}
\end{prop}

\begin{proof} When $d$ is even we have
\begin{align*}
T_\ast^\Sy(\ola\bfs)
=&\, \sum_{i=0}^d  \Big(\prod_{\ell=i+1}^{d} (-1)^{s_\ell} \Big) T_\ast(s_{i+1},\dots,s_d)T_\ast(s_i,\dots,s_1)\\
=&\, (-1)^w \sum_{i=0}^d  \Big(\prod_{\ell=1}^i (-1)^{s_\ell} \Big) T_\ast(s_i,\dots,s_1)T_\ast(s_{i+1},\dots,s_d)
=(-1)^w   T_\ast^\Sy(\bfs).
\end{align*}
The same argument works for SMSVs and odd $d$ cases.
\end{proof}

\begin{prop}
For all $s\in\N$
\begin{equation}\label{equ:Dblt}
t_\calA(s,s)=\zeta^{(2)}_\calA(s,s)=
\left\{
  \begin{array}{ll}
     \emtq_2^2/2,  &\quad \hbox{if $s=1$;} \\
    (1-2^{1-s})^2 \gb_s^2/2,  \quad &\quad \hbox{if $s\ge 2$,}
  \end{array}
\right.
\end{equation}
and
\begin{equation}\label{equ:Triplet}
t_\calA(s,s,s)=-\zeta^{(2)}_\calA(s,s,s)=
\left\{
  \begin{array}{ll}
     \emtq_2^3/6+\gb_3/8 ,  &\quad \hbox{if $s=1$;} \\
    (1-2^{1-s})^3 \gb_s^3/6+\gb_{3s}/8,  \quad &\quad \hbox{if $s\ge 2$.}
  \end{array}
\right.
\end{equation}
More generally, for all $d\in\N$
\begin{equation*}
t_\calA(\{s\}^d)=-\zeta^{(2)}_\calA(\{s\}^d)\in \gd_{s,1}\emtq_2^s\Q+\sum_{k_1+\dots+k_\ell=d,\ell\in\N} \gb_{sk_1}\dots \gb_{sk_\ell }\Q.
\end{equation*}
Moreover, we may assume all $k_j$'s are odd.
\end{prop}

\begin{proof}
By the stuffle relations,
\begin{align*}
2t_\calA(s,s)=&\, t_\calA(s)^2-t_\calA(2s).
\end{align*}
Thus \eqref{equ:Dblt} follows from \eqref{equ:FESdepth1odd} and \eqref{equ:FMtSVReversalEvenDep}
by noticing that $\gb_w=0$ if $w$ is even. Similarly, by the stuffle relations,
\begin{align*}
6t_\calA(s,s,s)=&\, t_\calA(s)t_\calA(s,s)-t_\calA(2s,s)-t_\calA(s,2s) \\
=&\,t_\calA(s)(t_\calA(s)^2s-t_\calA(2s))-t_\calA(s)t_\calA(2s)+t_\calA(3).
\end{align*}
Hence \eqref{equ:Triplet} follows from \eqref{equ:FESdepth1odd} and \eqref{equ:FMtSVReversalEvenDep}.
The general homogeneous FMtVs can be computed similarly by induction or by \cite[Theorem~2.3]{Hoffman2015}
(see also \cite[Lemma 5.1]{TaurasoZh2010}).
The statement for $\zeta^{(2)}_\calA$ is an easy application of the reversal relation \eqref{equ:FMtSVReversalEvenDep} so that the sign in the equation is $(-1)^{sd}$. But the values could be nonzero only when both $s$ and $d$ are odd, resulting in the negative sign.
\end{proof}

\subsection{Linear shuffle relations}
The most nontrivial relations among finite MZVs and finite Euler sums is provided by the linear shuffle relations
which is closely related to the shuffle relations among the classical MZVs and Euler sums. In this subsection,
we will extend this to FMMVs and their alternating versions.

For any $n\in\N$, $\bfs=(s_1,\dots,s_d)\in\N^d$,
$\bfeps=(\eps_1, \dots, \eps_r)\in\{\pm 1\}^r$, and $\bfgs=(\sigma_1, \dots, \sigma_r)\in\{\pm 1\}^r$ , we defined the
the partial sums of MMVs by
\begin{equation}\label{equ:AMMVdefn}
M_n(\bfs;\bfeps;\bfgs):=\sum_{n>k_1>\cdots>k_d>0} \frac{(1+\eps_1(-1)^{k_1})\sigma_1^{(2k_1+1-\eps_1)/4} \cdots (1+\eps_d(-1)^{k_d})\sigma_d^{(2k_d+1-\eps_d)/4}}{k_1^{s_1} \cdots k_d^{s_d}}.
\end{equation}
When $n\to \infty$ and $(s_d,\sigma_d)\neq (1,1)$ we recover the alternating MMVs studied by Xu, Yan and the last author
in \cite{XuZhao2020c,XuYanZhao2022Aug}. By taking $n$ to be primes, we can now define the \emph{finite alternating multiple mixed values}
\begin{equation}\label{equ:FAMMVdefn}
M_\calA(\bfs;\bfeps;\bfgs):=\Big( M_p(\bfs;\bfeps;\bfgs) \pmod{p} \Big)_{p\in\calP}.
\end{equation}

It turns out that in depth two, all the (non-alternating) MMVs have been given special names already,
which we now recall. Let $d=2$. Then we set
\begin{alignat*}{4}
\zeta^{(2)}_\calA(\bfs;\bfgs):=&\, M_\calA(\bfs;1,1;\bfgs), \qquad &
T_\calA(\bfs;\bfgs):=&\, M_\calA(\bfs;1,-1;\bfgs),\\
S_\calA(\bfs;\bfgs):=&\, M_\calA(\bfs;-1,1;\bfgs), \qquad  &
t_\calA(\bfs;\bfgs):=&\, M_\calA(\bfs;-1,-1;\bfgs).
\end{alignat*}
Moreover, to save space, if an alternating sign $\gs_j=-1$ then we put a bar on top of $s_j$ correspondingly.
For example,
\begin{equation*}
 S_\calA(\bar2,3) \sum_{m>n>0, 2\nmid m, 2|n} \frac{(-1)^{(m+1)/2}}{m^2n^3}.
\end{equation*}

We recall briefly the main setup for the integral expression of alternating MMVs. Let
\begin{equation*}
\mya=\frac{dt}{t}, \myb=w_{+1}^{-1}:=\frac{dt}{1-t^2},
\quad \myc=w_{-1}^{-1}:=\frac{-dt}{1+t^2},
\quad \gb=w_{+1}^{+1}:=\frac{t\, dt}{1-t^2},
\quad \gc=w_{-1}^{+1}:=\frac{-t\, dt}{1+t^2}.
\end{equation*}
Set
$$ w_{\gs}^{\eps_1,\eps_2}:=\max\{\gs,\sgn(1+\eps_2-\eps_1)\}w_{\gs}^{\eps_1\eps_2}. $$
Namely, $w_{\gs}^{\eps_1,\eps_2}=w_{\gs}^{\eps_1\eps_2}$ unless $\gs=\eps_2=-\eps_1=-1$
when $w_{\gs}^{\eps_1,\eps_2}=-w_{\gs}^{\eps_1\eps_2}$.
It is straight-forward to deduce that alternating MMVs can be expressed by the following iterated integrals
\begin{equation}\label{AMMV-iterated-integrals}
M(\bfs;\bfeps;\bfsi)=
\int_0^1 w_0^{s_1-1}w_{\gs_1}^{\eps_1,\eps_2}
w_0^{s_2-1}w_{\gs_1\gs_2}^{\eps_2,\eps_3}
\cdots
w_0^{s_r-1}w_{\gs_1\gs_2\cdots\gs_r}^{\eps_r}.
\end{equation}
For $\bfgs,\bfeps\in\{\pm 1\}^r$ define
\begin{equation*}
\sgn(\bfgs,\bfeps):=(-1)^{\sharp\{i<r \mid \gs_{i}=\eps_{i}=\eps_{i+1}\eps_{i+2}\cdots\eps_{r}=-1\}}.
\end{equation*}
Then for all $\bfs=(s_1,\dots,s_r)\in\N^r$ with $(s_1,\gs_1)\ne(1,1)$, we have
\begin{align}\label{AMMV-iterated-integrals-duality}
\int_0^1 w_0^{s_1-1}w_{\gs_1}^{\eps_1}\cdots w_0^{s_r-1}w_{\gs_r}^{\eps_r}
=\sgn(\bfgs,\bfeps)M(\bfs;\widetilde{\bfeps};\widetilde{\bfgs})
\end{align}
where $\widetilde{\bfgs}=(\gs_1, \gs_2\gs_1,\ldots, \gs_r\gs_{r-1})$ and
$\widetilde{\bfeps}=(\eps_1\cdots\eps_r, \eps_2\cdots\eps_r ,\ldots,\eps_{r-1}\eps_r, \eps_r)$.
See \cite{XuYanZhao2022Aug} for more details, where the definition of alternating MMVs differs from the one
used in this paper by a power of 2.

We now apply the above integral expressions to our finite situations. For example,
\begin{align*}
  \int_0^t \myb \gb \gc=&\,  \int_0^t  \sum_{i\ge0}  x^{2i}\, dx \int_0^x \sum_{j\ge0}  y^{2j+1}\, dy  \int_0^y \sum_{k\ge0} (-1)^{k+1} z^{2k+1}\, dz\\
=&\, \sum_{i,j,k\ge0} \frac{t^{2i+2j+2k+5} (-1)^{k+1}}{ (2i+2j+2k+5)(2j+2k+4)(2k+2)}     \\
=&\, \sum_{m>n>l>0, 2\nmid m, 2|n,2|l} \frac{t^{m} (-1)^{l/2}}{mnl}.
\end{align*}
Taking the coefficient of $t^p$ for any prime $p$, we get
\begin{align*}
 {\rm Coeff}_{t^p} \int_0^t \myb \gb \gc=  \frac{1}{p} \zeta^{(2)}_p(1,\bar1).
\end{align*}
Then from the shuffle relation
\begin{align*}
 \int_0^t \myb \int_0^t \gb \gc= \int_0^t \myb \sha \gb \gc =\int_0^t (\myb \gb \gc+\gb\myb \gc+\gb \gc\myb )
\end{align*}
we get
\begin{align*}
p\cdot \Coeff_{t^p} \Big(\int_0^t \myb \int_0^t \gb \gc \Big)
=\zeta^{(2)}_p(1,\bar1)+S_p(1,\bar1)+ t_p(\bar1,\bar1).
\end{align*}
Hence we arrive at one of the linear shuffle relations:
\begin{align} \label{equ:linSha1}
\zeta^{(2)}_\calA(1,\bar1)+S_\calA(1,\bar1)+t_\calA(\bar1,\bar1)=0.
\end{align}
Similarly, we can obtain the following linear shuffle relations by \eqref{AMMV-iterated-integrals-duality}:
\begin{align} \label{equ:linSha2}
\myb\sha\myb\myb &\, \Longrightarrow  3T_\calA(1,1)=0   ,\\
\myb\sha\myb\myc &\, \Longrightarrow  2T_\calA(1,\bar1)-T_\calA(\bar1,\bar1)=0,\label{equ:linSha3}\\
%\myb\sha\myc\myb &\, \Longrightarrow  T_\calA(\bar1,\bar1)+2(-1)^{p'}T_\calA(\bar1,1) =0 ,\label{equ:linSha4}\\
%\gb\sha\myb\gb &\, \Longrightarrow   S_\calA(1,1)+2\zeta^{(2)}_\calA(1,1)=0,\label{equ:linSha5}\\
%\gb\sha\gb \myb&\, \Longrightarrow   2t_\calA(1,1)+S_\calA(1,1)=0,\label{equ:linSha6}\\
\gc\sha\myb\gc &\, \Longrightarrow  (-1)^{p'}S_\calA(\bar1,\bar1)-2\zeta^{(2)}_\calA(\bar1,1)=0,\label{equ:linSha7}\\
\gc\sha\gc\myb &\, \Longrightarrow   2 t_\calA(1,\bar1) +S_\calA(\bar1,\bar1) =0,\label{equ:linSha8}\\
\myb\sha\gc\gb &\, \Longrightarrow  \zeta^{(2)}_\calA(\bar1,\bar1)-(-1)^{p'}S_\calA(\bar1,1)-(-1)^{p'} t_\calA(\bar1,1)=0. \label{equ:linSha9}
\end{align}
Here and in the rest of this section, we put $p'=(p-1)/2$ to save space.
We observe that the number of $\myb$ and $\myc$ must be either one or three in order to have nontrivial relations.

\section{Depth one and two values}
First we observe that since $\zeta_\calA(s)=0$ for all $s\in\N$, by \cite[Theorem~8.2.7]{ZhaoBook},
\begin{align}\label{equ:FESdepth1odd}
\zeta^{(2)}_\calA(s)=S_\calA(s)=-t_\calA(s)=-T_\calA(s)=\frac12\zeta_\calA(\bar{s})=&\,
\left\{
 \begin{array}{ll}
  -\tq_2,     \quad  & \quad \hbox{if $s=1$;}\\
  (2^{1-s}-1 )\gb_{s},\quad &\quad  \hbox{if $s\ge 2$,}
 \end{array}
\right.
\end{align}
where $\tq_2$ is the Fermat quotient \eqref{equ:FermatQ2} and $\gb_{s}$ is given by \eqref{equ:betaw}.
In weight one and depth one, we have the following result.

\begin{prop} \label{prop:D1wt1Value}
We have
\begin{alignat}{4}\label{equ:D1V1E1}
S_\calA(1)=&\,\zeta^{(2)}_\calA(1)=-\emtq_2, \quad &\, \zeta^{(2)}_\calA(\bar1)=&\,S_\calA(\bar1)=-\emtq_2/2, \\
t_\calA(1)=&\, T_\calA(1)=\emtq_2, \quad  &\, t_\calA(\bar1)=&\, T_\calA(\bar1)=-(-1)^{p'}\emtq_2/2 \label{equ:D1V1E2}
\end{alignat}
where we regard $(-1)^{p'}$ as $\big((-1)^{p'} \big)_{p\ge 3}\in\calA$.
\end{prop}

\begin{proof}
We only need to prove that $S_\calA(1)=2S_\calA(\bar1)$. Indeed,
\begin{alignat*}{4}
2S_p(1)=&\, \sum_{k=1}^{p'} \frac{(-1)^k+1}{k} + \sum_{p'<k<p,2|k} \frac{2}{k} \\
= &\, 2S_p(\bar1)+\sum_{k=1}^{p'} \frac{1}{k} + \sum_{0<k\le p',2\nmid k} \frac{1}{p-k}+ \sum_{p'<k<p,2|k} \frac{1}{k} \\
\equiv &\, 2S_p(\bar1)+\sum_{0<k\le  p',2|k} \frac{1}{k} + \sum_{p'<k<p,2|k} \frac{1}{k}  \quad & &  \pmod{p}\\
\equiv &\, 2S_p(\bar1)+S_p(1) \quad & &  \pmod{p}.
\end{alignat*}
The proposition follows easily from \eqref{equ:FESdepth1odd} and substitutions
$k\to p-k$ for $t_\calA(\bar1)$ and $T_\calA(\bar1)$.
\end{proof}

For depth two values, we have the following results.

\begin{prop} \label{prop:Depth2Value}
For all $a,b\in\N$, if $w=a+b$ is odd then
\begin{align}\label{equ:DblFMtVs}
\zeta^{(2),\star}_\calA(a,b)=&\,-t^\star_\calA(a,b)=  \frac12\Big[2^{1-w}-1-(-1)^a 2^{-w}\binom{w}{a}\Big]\gb_w, \\
\zeta^{(2)}_\calA(a,b)=&\, -t_\calA(a,b)=\frac12\Big[1-2^{1-w}-(-1)^a 2^{-w}\binom{w}{a}\Big]\gb_w, \\
S_\calA(a,b)=&\,T_\calA(a,b)=\frac{(-1)^a} 2  \Big(1-2^{-w}\Big) \binom{w}{a} \gb_w.\label{equ:DblFMSTVs}
\end{align}
\end{prop}

\begin{rem}
The formulas for $\zeta^{(2)}_\calA(s)$ in \eqref{equ:FESdepth1odd} and for $\zeta^{(2)}_\calA(a,b)$ in \eqref{equ:DblFMtVs}
are consistent with \cite[Proposition~2.1]{KanekoMuYo2021}. Note that the ordering in this paper is opposite to that of \cite{KanekoMuYo2021}.
Also, our definition of $\zeta^{(2)}_\calA(\bfs)$ is $2^{-|\bfs|}$ times that in \cite{KanekoMuYo2021}.
\end{rem}

\begin{proof}
By \cite[Theorem~8.6.4]{ZhaoBook},
\begin{align} \notag
\zeta^\star_\calA(a,b)=&\,\zeta_\calA(a,b)= (-1)^a\binom{w}{a}\gb_w,  \\
\zeta^\star_\calA(\ol{a},\ol{b})=&\,\zeta_\calA(\ol{a},\ol{b})= (-1)^a(2^{1-w}-1)\binom{w}{a} \gb_w, \notag \\
\zeta_\calA(\ol{a},b)=&\, \zeta_\calA(a,\ol{b})= (1-2^{1-w})\gb_w,  \label{equ:zetaBarab}\\
\zeta^\star_\calA(\ol{a},b)=&\, \zeta^\star_\calA(a,\ol{b})= (2^{1-w}-1)\gb_w.\notag
\end{align}
The proposition follows easily.
\end{proof}

\section{Weigh two finite alternating MMVs}
We first recall that the Euler polynomials $E_n(x)$ are defined by the generating function
\begin{equation*}
\frac{2 e^{tx} }{e^t+1} =\sum_{n=0}^\infty E_n(x) \frac{t^n}{n!}.
\end{equation*}
Moreover, $E_0(0)=1$ and for all $j\in \N$
\begin{equation}\label{equ:EulerBern}
E_j(0) = \frac{2^{j+1}}{j+1} \Big(B_{j+1}\Big(\frac{1}2\Big)
-B_{j+1}\Big)= \frac{2}{j+1}(1-2^{j+1})B_{j+1}
\end{equation}
by \cite[p.~242, (8.11)]{ZhaoBook}. The Euler numbers $E_n$ are  defined by the generating function
\begin{equation*}
\frac{2}{e^t+e^{-t}} =\sum_{n=0}^\infty E_n \frac{t^n}{n!},
\end{equation*}
which then satisfy
\begin{equation}\label{equ:EulerNumber}
E_n =  2^n  E_n\Big(\frac12 \Big).
\end{equation}
Let $G$ be the traditional Catalan's constant. Then
\begin{equation*}
   T(\bar2)=\sum_{n\ge 1} \frac{(-1)^{n}}{(2n-1)^2} =-G.
\end{equation*}
Motivated by the next proposition, we define the \emph{finite Catalan's constant} by
\begin{equation*}
     \tG_\calA:=\bigg( \frac{E_{p-3}}2 \bigg)_{3<p\in\calP} \in\calA.
\end{equation*}

\begin{prop}\label{prop:Tbar2}
We have
\begin{equation}\label{equ:Tbar2}
T_\calA(\bar2)= - \tG_\calA.
\end{equation}
\end{prop}
\begin{proof} For any prime $p\ge 5$, we have modulo $p$
\begin{align*}
T_p(\bar2)=&\, \sum_{0<k<p, 2\nmid k} \frac{(-1)^{(k+1)/2}}{k^2}
=\sum_{0<k<p, 2| k} \frac{(-1)^{(p-k+1)/2}}{(p-k)^2} \\
\equiv &\, (-1)^{(p+1)/2} \sum_{0<k<p, 2| k} \frac{(-1)^{k/2}}{k^2}   \\
\equiv &\, (-1)^{(p+1)/2} \sum_{n=1}^{p'} \frac{(-1)^n}{4 n^2}   \\
\equiv &\,\frac{(-1)^{(p+1)/2}}{4}  \sum_{n=1}^{p'}  (-1)^n n^{p-3}  \\
\equiv &\, \frac{(-1)^{(p+1)/2}}{8} \Big((-1)^{p'} E_{p-3}\Big(\frac{p+1}{2} \Big)+E_{p-3}\big(0\big)\Big)   \\
\equiv &\, -\frac{1}{8}   E_{p-3}\Big(\frac{1}{2} \Big)  \\
\equiv &\, -\frac{1}{2}   E_{p-3}
\end{align*}
by \cite[Lemma 8.2.5]{ZhaoBook} and \eqref{equ:EulerNumber}
since $E_{p-3}(0)=0$ for all odd primes $p$ by \eqref{equ:EulerBern}.
This completes the proof of the proposition.
\end{proof}

\begin{prop}\label{prop:wt2FMTVs}
We have
\begin{equation}\label{equ:wt2FMTVs}
S_\calA(2)=T_\calA(2)=T_\calA(1,1)=0, \quad S_\calA(1,1)= - \emtq_2^2, \quad t_\calA(1,1)=\zeta^{(2)}_\calA(1,1)= \frac{\emtq_2^2}2.
\end{equation}
\end{prop}
\begin{proof}
The vanishing of the first three values follows from \eqref{equ:FESdepth1odd} and \eqref{equ:linSha2} quickly.
Then by Proposition~\ref{prop:D1wt1Value} we get
\begin{equation*}
S_\calA(1,1)=S_\calA(1,1)+T_\calA(1,1) =S_\calA(1)T_\calA(1)=\tq_2^2.
\end{equation*}
By \eqref{equ:FESdepth1odd} we see that
\begin{equation*}
  2 t_\calA(1,1)=t_\calA(1)^2-t_\calA(2)= \tq_2^2
\end{equation*}
and similarly,
\begin{equation*}
  2 \zeta^{(2)}_\calA(1,1)=S_\calA(1)^2-S_\calA(2)= \tq_2^2.
\end{equation*}
We have completed the proof of the proposition.
\end{proof}

\begin{prop}\label{prop:wt2FMTVbars}
We have
\begin{equation}\label{equ:wt2FMTVbars}
T_\calA(\bar1,\bar1)=2T_\calA(1,\bar1), \quad S_\calA(\bar1,\bar1)= -2 t_\calA(1,\bar1), \quad
t_\calA(\bar1,\bar1)=\zeta^{(2)}_\calA(\bar1,\bar1)= \frac{\emtq_2^2}8.
\end{equation}
\end{prop}
\begin{proof}
The first identity is just \eqref{equ:linSha3} and the second is \eqref{equ:linSha8}.
The other two follows immediately from the stuffle relations by Proposition~\ref{prop:D1wt1Value}.
\end{proof}

\begin{prop}\label{prop:wt2FMTVall}
We have
\begin{alignat}{4}\label{equ:wt2FMTVall}
T_\calA(\bar2)=&\, -\tG_\calA,
\quad T_\calA(\bar1,1)= -\frac{(-1)^{p'}} 2 \tG_\calA,   &\,
\quad S_\calA(\bar1,1)=&\, \frac{(-1)^{p'}}2 \emtq_2^2-\frac12 \tG_\calA,\\
S_\calA(\bar2)=&\, (-1)^{p'} \tG_\calA,
\quad T_\calA(1,\bar1)=\frac12 \tG_\calA, &\,
\quad S_\calA(1,\bar1)=&\,-\frac{1}2\emtq_2^2+\frac{(-1)^{p'}}2 \tG_\calA ,\\
\quad t_\calA(\bar1,1)=&\, -\frac{3(-1)^{p'}}8 \emtq_2^2+\frac12 \tG_\calA, &\,
\quad \zeta^{(2)}_\calA(\bar1,1)=&\, \frac18\emtq_2^2 -\frac{(-1)^{p'}}2\tG_\calA,\\
\quad t_\calA(1,\bar1)=&\, -\frac{(-1)^{p'}}8 \emtq_2^2 + \frac12 \tG_\calA, &\,
\quad \zeta^{(2)}_\calA(1,\bar1)=&\, \frac{3}8\emtq_2^2-\frac{(-1)^{p'}}2 \tG_\calA.
\end{alignat}
\end{prop}

\begin{proof}
Fix a large prime $p$ so that identities in both Proposition~\ref{prop:wt2FMTVs} and Proposition~\ref{prop:wt2FMTVbars}
hold for $T_p$, $S_p$, $t_p$, and $\zeta^{(2)}_p$.
Throughout the rest of this proof, we will drop the subscript $p$ to save space.

By stuffle relation and Proposition~\ref{prop:D1wt1Value} we have
\begin{align}\label{equ:-1suppq2}
(-1)^{p'}\frac{\tq_2}2=&\, S(1)T(\bar1)=S(\bar1,1) + T(1,\bar1),  \\
(-1)^{p'}\frac{\tq_2}4=&\, S(\bar1)T(\bar1)=S(\bar1,\bar1) + T(\bar1,\bar1)
=2T(1,\bar1)-2 t(1,\bar1), \label{equ:-1suppq}
\end{align}
by Proposition~\ref{prop:wt2FMTVbars}. Plugging \eqref{equ:-1suppq} into \eqref{equ:-1suppq2},
\begin{equation}\label{equ:Sbar11}
S(\bar1,1)=3T(1,\bar1)-4t(1,\bar1).
\end{equation}
By changing of indices $k\to p-k$ for $\zeta^{(2)}$ in \eqref{equ:linSha9} (or by \eqref{equ:linSha1}), we get
\begin{equation}\label{equ:tbar11}
t(\bar1,1) =(-1)^{p'}t(\bar1,\bar1)-S(\bar1,1)
=(-1)^{p'}\frac{\tq_2^2}8-S(\bar1,1)
=3t(1,\bar1)-2T(1,\bar1)
\end{equation}
by \eqref{equ:-1suppq} and \eqref{equ:Sbar11}.

Next, by Proposition~\ref{prop:D1wt1Value} and \eqref{equ:-1suppq}
\begin{equation*}
T(\bar2)=t(1)t(\bar1)-t(1,\bar1)-t(\bar1,1)=-\frac{(-1)^{p'}}2 \tq_2^2-t(1,\bar1)-t(\bar1,1)
=-2T(1,\bar1)
\end{equation*}
by \eqref{equ:tbar11} and \eqref{equ:tbar11}. Thus
\begin{equation*}
T_\calA(1,\bar1)=-\frac12 T_\calA(\bar2)= \frac12 \tG_\calA
\end{equation*}
by Proposition~\ref{prop:Tbar2}. On the other hand, by changing index $k\to p-k$
we see that
\begin{equation*}
S_\calA(\bar2)=-(-1)^{p'}  T(\bar2)=(-1)^{p'}\tG_\calA,\qquad
T_\calA(\bar1,1)=-(-1)^{p'}T_\calA(\bar1,1)= -\frac{(-1)^{p'}}2\tG_\calA.
\end{equation*}

Further, by \eqref{equ:linSha7} we see that
\begin{equation}\label{equ:zetabar11}
\zeta^{(2)}(\bar1,1)=\frac{(-1)^{p'}}2 S(\bar1,\bar1)
= \frac{\tq_2^2}8- \frac{(-1)^{p'}}2 T(\bar1,\bar1)
= \frac{\tq_2^2}8-(-1)^{p'}T(1,\bar1)
= \frac{\tq_2^2}8 +T(\bar1,1)
\end{equation}
by \eqref{equ:-1suppq} and Proposition~\ref{prop:wt2FMTVbars}. Hence
\begin{equation}\label{equ:t1bar1}
t(1,\bar1)=-(-1)^{p'}\zeta^{(2)}(\bar1,1)
= -\frac{(-1)^{p'}}8 \tq_2^2 + T(1,\bar1)= -\frac{(-1)^{p'}}8\tq_2^2 + \frac12 \tG_\calA .
\end{equation}
All the other identities in the proposition can be now derived by applying the change of indices $k\to p-k$.
\end{proof}

By combining Propositions \ref{prop:wt2FMTVs}, \ref{prop:wt2FMTVbars} and \ref{prop:wt2FMTVall}
we immediately obtain the following theorem.
\begin{thm}
Let $\FAMMV_w$ be the $\Q$-vector space generated by finite alternating MMVs of weight $w$.
Let $(-1)^{p'}=\big((-1)^{(p-1)/2} \big)_{3\le p\in\calP}\in\calA$.
Then
\begin{equation}\label{equ:FAMbasis}
\FAMMV_2=\langle \emtq_2^2, \tG_\calA^2, (-1)^{p'}\emtq_2^2, (-1)^{p'}\tG_\calA^2 \rangle.
\end{equation}
\end{thm}

We have carried out some extensive computations of $\FAMMV_w$ for small weights $w$
and tabulate the result in the Appendix to this paper.

\section{Sum formulas of symmetric and finite MTVs of even depth}
For all $w\ge d\in\N$ let
\begin{equation*}
I_{w,d}:= \{\bfs: |\bfs|=w,\dep(\bfs)=d\}.
\end{equation*}

\begin{prop}
Suppose $w\ge d\in\N$ with $w$ odd and $d$ even. Then for $F=T$ and $S$ we have
\begin{equation*}
\sum_{\bfs\in I_{w,d}} F_\calA(\bfs)=0, \quad \sum_{\bfs\in I_{w,d}} F_\ast^\Sy(\bfs)=0
\end{equation*}
\end{prop}

\begin{proof}
By reversal relation we have
\begin{equation*}
\sum_{\bfs\in I_{w,d}} T_\calA(\bfs)=\sum_{\bfs\in I_{w,d}} (-1)^w T_\calA(\ola{\bfs})=-\sum_{\bfs\in I_{w,d}} T_\calA(\bfs)=0.
\end{equation*}
The same argument works for FMSVs. For the symmetric values, by \eqref{equ:SMTSVReversalEvenDep} we have
\begin{align*}
\sum_{\bfs\in I_{w,d}} T_\ast^\Sy(\bfs)=\sum_{\bfs\in I_{w,d}} T_\ast^\Sy(\ola\bfs)
=(-1)^w\sum_{\bfs\in I_{w,d}} T_\ast^\Sy(\bfs)=-\sum_{\bfs\in I_{w,d}} T_\ast^\Sy(\bfs)=0.
\end{align*}
The same argument works for SMSVs. This concludes the proof of the proposition.
\end{proof}

\section{Restricted sum relations}

\begin{thm}
Let $w,d\in\N$ with $d\le w$. Let
\begin{align*}
I_{w,d,i}:=&\, \{\bfs: |\bfs|=w,\dep(\bfs)=d,s_i\ge 2\} \forall 1\le i\le d,\\
I_{w,d,i,i}:=&\, \{\bfs: |\bfs|=w,\dep(\bfs)=d,s_i\ge 3\} \forall 1\le i\le d, \\
I_{w,d,i,j}:=&\, \{\bfs: |\bfs|=w,\dep(\bfs)=d,s_i\ge 2,s_j\ge 2\} \forall 1\le j<i\le d.
\end{align*}
Then
\begin{align*}
    \sum_{\bfs\in I_{w,d}} M_\calA(\bfs)=&\,\sum_{\bfs\in I_{w,d}} M^\star_\calA(\bfs)= 0,\\
    \sum_{\bfs\in I_{w,d,i}} M_\calA(\bfs)=&\, (-1)^{i-1}\left(\binom{w-1}{i-1}+(-1)^d\binom{w-1}{d-i}\right)\gb_w,\\
    \sum_{\bfs\in I_{w,d,i}} M^\star_\calA(\bfs)=&\, (-1)^{i-1}\left((-1)^d\binom{w-1}{i-1}+\binom{w-1}{d-i}\right)\gb_w,\\
 \sum_{\bfs\in I_{w,d,i,j}} M_\calA(\bfs)=&\,  (-1)^d\sum_{\bfs\in I_{w,d,i,j}} M^\star_\calA(\bfs)= \frac12N_{w,d,i,j}\gb_w,
\end{align*}
where $\gb_w$ is defined by \eqref{equ:betaw} and
$N_{w,d,i,j}$ is an integer explicitly given by \cite[Theorem~3.1]{MuraharaSa2019}.
\end{thm}

\begin{proof}
For any prime $p$ and $\bfs\in I_{w,d}$, by \cite[Theorem~4.1]{Hoffman2015} we have
\begin{equation*}
\sum_{\gs\in \frakS_d} H_p(\gs(\bfs))\equiv \sum_{\gs\in \frakS_d} H^\star_p(\gs(\bfs))\equiv 0 \pmod{p}
\end{equation*}
where $\frakS_d$ is the group of symmetry of $d$ letters. By partitioning $I_{w,d}$ into
equivalent classes under permutation we see quickly that
\begin{equation*}
\sum_{\bfs\in I_{w,d}} \zeta_\calA(\bfs)=0.
\end{equation*}
Thus,
\begin{equation*}
\sum_{\bfs\in I_{w,d}} M_\calA(\bfs)
=\sum_{\bfs\in I_{w,d}} \sum_{0<n_1<\cdots<n_d<p} \prod_{j=1}^d \frac{(1+(-1)^{n_j})+(1-(-1)^{n_j})}{2n_j^{s_j}}
= \sum_{\bfs\in I_{w,d}} \zeta_\calA(\bfs)=0.
\end{equation*}
The second and third equations follow from the sum formulas of Saito and Wakabayashi \cite[Theorem~1.4]{SaitoWa2013b}
and the last two from \cite[Theorem~3.1]{MuraharaSa2019} immediately.
\end{proof}

\bigskip

\begin{center}
{\Large \bf Appendix: Dimensions of subspaces of finite MMVs} 
\end{center}
\begin{center}
by {\sc Jeremy Feng, Angelina Kim, Sienna Li, Ryan Qin,  Logan Wang, and J.\ Zhao}
\end{center}

In the follow table we provide the conjectural dimension of various subspaces of $\FMMV$ considered in this paper.
We achieve this by numerically finding all the $\Q$-linear relations
among these values aided by Maple computation using the code contained in \cite[Appendix D]{ZhaoBook}. 

\begin{table}[!h]
{
\begin{center}
\begin{tabular}{  |c|c|c|c|c|c|c|c|c|c|c|c|c|c|c| } \hline
       $w$     & 0 &  1  &  2  &  3  &  4  &  5  &  6  & 7  &  8  &  9  &  10 &  11  &  12  &  13  \\ \hline
$\dim_\Q \FMZV_w$  & 1 &  0  & 0  &  1  &  0  &  1  & 1 &  1 &  2 &  2 &  3  &  4 &  5&  7\\ \hline
$\dim_\Q\FMZV^{(2)}_w$ &0 & 1 &  1  &  2  &  3  &  5  &  8  & 13 &  21 &  34 &  55 &  89  &  144 &  233 \\ \hline
$\dim_\Q \FES_w$  &0 & 1 &  1  &  2  &  3  &  5  &  8  & 13 &  21 &  34 &  55 &  89  &  144 &  233 \\ \hline
$\dim_\Q \FMtV_w$  &0 & 1 &  1  &  2  &  3  &  5  &  8  & 13 &  21 &  34 &  55 &  89  &  144 &  233 \\ \hline
$\dim_\Q \FMTV_w$  & 0& 1 &  0  &  1  &  2  &  3  &  3  &  6  & 9  &  15  &  17 &  32 &  44 &  76  \\ \hline
$\dim_\Q \FMSV_w$  & 0 &  1  &  1  &  1  &  2  &  4  &  5  & 7 &  12 &  19 & 28  &  39  & 66   &  ?  \\ \hline
$\dim_\Q \FMMV_w$  & 0 &  1  &  1  &  2  &  3  &  5  &  8  & 13 &  21 &  34 &  55  &  89 ? &  144 ?&  233? \\ \hline
$\dim_\Q \FMMVe_w$  & 0 &  1  &  1  &  2  &  3  &  5  &  8  & 13 &  21 &  34 &  55 &  89 ? &  144? &  233? \\ \hline
$\dim_\Q \FMMVo_w$ & 0 &  1  &  1  &  2  &  3  &  5  &  8  & 13 &  21 &  34 &  55  &  89 ? &  144? &  233? \\ \hline 
\end{tabular}
\end{center}
}
\caption{Conjectural Dimensions of Various Subspaces of $\MMV$.}
\label{Table:dimMMV}
\end{table}

 From numerical computation we can formulate the following conjecture.
\begin{conj}\label{conj:dimFMMVs}
(i) For all $w\ge 1$,
\begin{align*}
\FMZV^{(2)}_w =\FMMVe_w=\FMMV_w=\FES_w=\FMMVo_w=\FMtV_w
\end{align*}
all have dimension $F_w$.

(ii) For all $w\ge 1$,
\begin{align*}
\frac{\FMZV^{(2)}_w}{\zeta^{(2)}_\calA(\{1\}^w)\Q}
=\frac{\FMMVe_w}{\zeta^{(2)}_\calA(\{1\}^w)\Q}
=\frac{\FES_w}{\zeta_\calA(\{\bar1\}^w)\Q}
=\frac{\FMMVo_w}{t_\calA(\{1\}^w)\Q}
=\frac{\FMtV_w}{t_\calA(\{1\}^w)\Q}
\end{align*}
all have dimension $F_w-1$.

(iii) For all $k\ge 1$,
\begin{align*}
\dim_\Q \FMTV_{2k+1}=\dim_\Q \FMTV_{2k}+\dim_\Q \FMTV_{2k-1}.
\end{align*}
\end{conj}

Note that Conjecture \ref{conj:dimFMMVs}(ii) and \cite[Theorem~7.1]{XuZhao2020a}
can partially explain why the conjectured dimension for $\FMMV_w$ (and $\FMMVe_w$) computed in this paper differs by 1 from
the corresponding dimension for classical $\MMV_w$ (and $\MMVe_w$) numerically computed in \cite{XuZhao2020a}.

%\medskip
%{\bf Acknowledgments.}  Jianqiang Zhao is supported by the Jacobs Prize from The Bishop's School.

\end{document}